\documentclass[12pt]{article}

\usepackage{tabularx}
\usepackage{graphics}
\usepackage{rotating}
\usepackage{lscape}
\usepackage{amsfonts,amssymb,stmaryrd}

\usepackage{pstricks}
\usepackage{pst-node}

\newtheorem{prop}{Proposition}[section]

\newenvironment{proof}{\noindent {\sc Proof.}}{\hfill$\square$}

\newenvironment{definition}{\noindent {\bf Definition }}{}
\newenvironment{remark}{\noindent {\bf Remark }}{}

\def\longto{\longrightarrow}

\def\ZZ{\mathbb{Z}}
\def\NN{\mathbb{N}}
\def\RR{\mathbb{R}}
\def\QQ{\mathbb{Q}}
\def\Isom{{\rm Isom}}
\def\GL{{\rm GL}}
\def\GA{{\rm GA}}
\def\Aut{{\rm Aut}}
\def\Id{{\rm Id}}
\def\om{\omega}

\def\P{{\mathcal P}}
\def\D{{\mathcal D}}
\def\F{{\mathcal F}}
\def\Q{{\mathcal Q}}
\def\C{{\mathcal C}}
\def\CC{\mathcal{CC}}
\def\S{{\mathcal S}}

\def\L{\Lambda}
\def\hLambda{{\hat{\Lambda}}}
\def\hL{\hLambda}
\def\vf{{\overrightarrow{f}}}
\def\Sym{\Sigma}

\def\cprlp{ primitive centered regular lattice polytope }
\def\rlp{ regular lattice polytope  }

\begin{document}

\title{Regular Lattice Polytopes and Root Systems}
\author{Pierre-Louis Montagard and Nicolas Ressayre}

\maketitle
\section{Introduction}
Let $\L$ be a  lattice in a real finite dimensional vector space. 
Here, we are interested in the lattice polytopes, that is the convex hulls of finite
subsets of $\L$.
Consider the group $G$ of the affine real transformations which map $\L$ onto itself. 
 Replacing the group of euclidean motions by the group $G$ one can define the notion of regular lattice polytopes. 
More precisely, for a lattice polytope $\P$, we denote by  $\Isom(\P)$ the subgroup of $G$ which preserves $\P$ and 
$\P$ is said to be a {\it regular lattice polytope} if the group $\Isom(\P)$ acts transitively 
on the set of complete flags of $\P$.
In \cite{Kar}, Karpenkov obtained a classification of the regular lattice polytopes. 
Here we obtain this classification by a more conceptual method. 
Another difference is that Karpenkov uses in an essential way 
the classification of the euclidean regular polytopes, but we don't.

Let us explain our approach. 
Firstly, we associate in a very natural way a reduced simply laced root system 
(not necessarily irreducible) to any regular lattice polytope $\P$. 
Then considering the faces of $\P$, we even show that the only possible root 
systems are of type $A_n$, $D_n$, $E_6,E_7,E_8$ and $(A_1)^n$ 
(later, we show that the exceptional root systems do not occur).
Conversely, we fix such a root system $\Phi$ and seek all the regular lattice polytopes 
$\P$ with $\Phi$ as associated root system. 
Such a polytope is characterized up to isomorphism by a lattice between the root lattice and the weight 
lattice, and a dominant weight.
We obtain in this way the list presented in Table~1.
 
For the convenience of the reader we also present below the regular lattice 
polytopes of dimension two in Figures~\ref{fig:hexagone}, \ref{fig:triangle} and~\ref{fig:carre}.
In each figure, we have two lattices: the intersection of the gray lines and the marked points. 
These lattices are the weight lattices $\L_P$ and the root lattices $\L_R$ of of the root system $A_2$
in the two first case and $A_1\times A_1$ in the last one.
In Figure~\ref{fig:hexagone}, we have drawn an hexagon which can be considered as a lattice polygon in 
$\L_P$ or $\L_R$: this gives two classes of regular hexagons. 
In Figure~\ref{fig:carre}, the situation is similar with  squares instead hexagons.
In Figure~\ref{fig:triangle}, we have two triangles: the dashed one in $\L_P$ and the other one in $\L_R$.
The result in dimension 2 asserts  that up to evident equivalence (see Section~\ref{sec:eq}) the only regular lattice  
polygons are these 2 hexagons, these 2 triangles and these 2 squares.

\psset{unit=.6,linecolor=lightgray,dotscale=3,dotsize=4pt,dotstyle=*}

\begin{figure}[!h]
\begin{tabularx}{\textwidth}{XXX}

  \begin{pspicture}(-4,-3)(4,3)
\pspolygon(-3,0)(-1.5,2.598)(1.5,2.598)(-1.5,-2.598)
\pspolygon(3,0)(1.5,-2.598)(-1.5,-2.598)(1.5,2.598)
\psline(1.5,2.598)(-1.5,2.598)
\psline(-3,0)(3,0)
\pspolygon(-2.5,-0.866)(-0.5,2.598)(2,-1.732)(-2,-1.732)(0.5,2.598)(2.5,-0.866)
\pspolygon(2.5,0.866)(0.5,-2.598)(-2,1.732)(2,1.732)(-0.5,-2.598)(-2.5,0.866)
\psline(-1.5,2.598)(1.5,-2.598)

\psdot[linecolor=black](0,0)
\psdot[linecolor=black](1.5,0.866)\psdot[linecolor=black](-1.5,0.866)\psdot[linecolor=black](1.5,-0.866)\psdot[linecolor=black](-1.5,-0.866)
\psdot[linecolor=black](0,1.732)\psdot[linecolor=black](0,-1.732)
\psdot[linecolor=black](3,0)\psdot[linecolor=black](-3,0)
\psdot[linecolor=black](-1.5,2.598)\psdot[linecolor=black](-1.5,-2.598)\psdot[linecolor=black](1.5,-2.598)\psdot[linecolor=black](1.5,2.598)

\pspolygon[linecolor=black,linewidth=0.05](1.5,0.866)(0,1.732)(-1.5,0.866)(-1.5,-0.866)(0,-1.732)(1.5,-0.866)
  \end{pspicture}
\caption{Two Hexagons}
  \label{fig:hexagone}
& 
\begin{pspicture}(-4,-3)(4,3)
\pspolygon(-3,0)(-1.5,2.598)(1.5,2.598)(-1.5,-2.598)
\pspolygon(3,0)(1.5,-2.598)(-1.5,-2.598)(1.5,2.598)
\psline(1.5,2.598)(-1.5,2.598)
\psline(-3,0)(3,0)
\pspolygon(-2.5,-0.866)(-0.5,2.598)(2,-1.732)(-2,-1.732)(0.5,2.598)(2.5,-0.866)
\pspolygon(2.5,0.866)(0.5,-2.598)(-2,1.732)(2,1.732)(-0.5,-2.598)(-2.5,0.866)
\psline(-1.5,2.598)(1.5,-2.598)

\psdot[linecolor=black](0,0)
\psdot[linecolor=black](1.5,0.866)\psdot[linecolor=black](-1.5,0.866)\psdot[linecolor=black](1.5,-0.866)\psdot[linecolor=black](-1.5,-0.866)
\psdot[linecolor=black](0,1.732)\psdot[linecolor=black](0,-1.732)
\psdot[linecolor=black](3,0)\psdot[linecolor=black](-3,0)
\psdot[linecolor=black](-1.5,2.598)\psdot[linecolor=black](-1.5,-2.598)\psdot[linecolor=black](1.5,-2.598)\psdot[linecolor=black](1.5,2.598)

\pspolygon[linecolor=black,linewidth=0.05](3,0)(-1.5,2.598)(-1.5,-2.598)
\pspolygon[linecolor=black,linewidth=0.05,linestyle=dashed](1,0)(-0.5,0.866)(-0.5,-0.866)
  \end{pspicture}
  \caption{Two Triangles}
  \label{fig:triangle}
 &  
  \begin{pspicture}(-4,-3)(4,3)
\pspolygon(-2,-2)(-2,2)(2,2)(2,-2)
\psline(-1,-2)(-1,2)\psline(0,-2)(0,2)\psline(1,-2)(1,2)
\psline(-2,-1)(2,-1)\psline(-2,0)(2,0)\psline(-2,1)(2,1)

\psdot[linecolor=black](0,0)
\psdot[linecolor=black](1,1)\psdot[linecolor=black](-1,1)\psdot[linecolor=black](1,-1)\psdot[linecolor=black](-1,-1)
\psdot[linecolor=black](2,0)\psdot[linecolor=black](2,2)\psdot[linecolor=black](0,2)\psdot[linecolor=black](-2,0)\psdot[linecolor=black](-2,-2)\psdot[linecolor=black](0,-2)
\psdot[linecolor=black](-2,2)\psdot[linecolor=black](2,-2)

\pspolygon[linecolor=black,linewidth=0.05](1,1)(-1,1)(-1,-1)(1,-1)
  \end{pspicture}
  \caption{Two Squares}
  \label{fig:carre}
\end{tabularx}
\end{figure}

Finally we briefly mention the well-known link between convex polytopes and algebraic geometry. 
We do not use this link but our inspiration for some results are of geometric origin. 
To each lattice polytope $\P$ one can associate a toric variety $X_\P$ with $T$ as torus, see for example \cite{Oda}.
The polytope $\P$ is regular if and only if the group of regular toric automorphism of $X_\P$ acts transitively on
the set of the maximal chains of irreducible $T$-stable subvarieties of $X_\P$.
In \cite{Pro}, Procesi consider the toric variety $X_\Phi$ associated to the 
decomposition in Weyl chambers of the root system $\Phi$. 
Our first results(see Proposition~\ref{prop:Wtrans}) can be
translated in the following way: there exists an equivariant surjective 
morphism from $X_\Phi$ onto $X_\P$ if 
$\Phi$ is the root system associated to the regular lattice polytope $\P$.\\

\noindent
{\bf Convention}
In this paper, we only consider non degenerated polytopes that is which 
span affinely the ambient real vector space.

\section{An equivalence relation}
\label{sec:eq}

Let $\Lambda$ be a free abelian group of rank $n$.
Let $\hLambda$ be a set with a free transitive action of $\Lambda$ denoted
by: $m+z$ for any $m\in\hLambda$ and $z\in\Lambda$.
Such a set $\hLambda$ is called a {\it $\Lambda$-affine space}.
A map $f\,:\,\hLambda\longto\hLambda$ is said to be {\it affine} if there exists a group morphism $\vf\,:\,\Lambda\longto\Lambda$ such that 
$f(m+z)=f(m)+\vf(z)$.

Let $\GL(\Lambda)\simeq\GL_n(\ZZ)$ denote the automorphism group of 
$\Lambda$ and $\GA(\hLambda)$ be the group of bijective affine maps of 
$\hLambda$.
We have the following split exact sequence:

$$
1\:\longto \:\L\longto\:\GA(\hL)\:\longto\:\GL(\L)\:\longto\:  1.
$$

Consider $\Lambda_\RR:=\Lambda\otimes\RR$ and its affine space 
$\hLambda_\RR:=m+\Lambda_\RR$ (for any $m\in\hLambda$). Now, $\hLambda$ 
is a lattice in $\hLambda_\RR$ and $\GA(\hLambda)$ is the subgroup of the 
isomorphism group $\GA(\hLambda_\RR)$ of $\hLambda_\RR$ of the elements 
which map $\hLambda$ onto itself.

A lattice polytope $\P$ is the convex hull in $\hLambda_\RR$ of a finite subset 
of $\hLambda$. 
Set
$$
\Isom(\P)=\{g\in\GA(\hLambda)\ |\ g\P=\P\}.
$$
A lattice polytope $\P$ is said to be {\it regular} if $\Isom(\P)$ acts
transitively on the set of complete flags of $\P$. We want to classify 
the regular polytope modulo $\GA(\hLambda)$ of course, but there is another reduction. 
We note  that if $h$ is a homothety of center in $\hLambda$ and integer ratio then $h$ normalize 
$\GA(\hLambda)$ and if $\P$ is a regular lattice polytope so is $h(\P)$. 
Finally we want to classify the regular lattice polytopes up to the group generated by $\GA(\hLambda)$ 
and the homotheties of center in $\hLambda$ and integer ratio. So we define:\\

\begin{definition}
We call $\mathcal H$ the subgroup of $\GA(\hLambda_\RR)$ generated by $\GA(\hLambda)$ 
and the homotheties of center in $\hLambda$ and integer ratio.\\
\end{definition}

Actually, this group doesn't  acts on the set of lattice polytopes, but on those with rational vertices. 
Nevertheless, this group defines a equivalence relation on the lattice polytopes. 
Our first reduction is to choose a common origin for the polytopes. 
More precisely let us fix an origin $O$ in $\hLambda$; now, one can 
identify $\Lambda$ and $\hLambda$, and  embed $\GL(\Lambda)$ in $\GA(\hLambda)$). 
Now we can define:\\   

\begin{definition}\label{def:redcentprim}
We call a lattice polytope {\it centered} if its barycenter is $O$; 
it is said to be {\it primitive} if it is not the image of  another polytope by  
an homothety of center $O$ and of integer ratio bigger than one.\\ 
\end{definition}

The following proposition reduce the classification to those of primitive centered
regular lattice polytope.

\begin{prop}
\label{prop:redcentprim}
  \begin{enumerate}
  \item 
  Let $\P$ be a lattice polytope. There exits $g\in\mathcal H$ such that $g(\P)$ is 
a primitive centered regular lattice polytope.
  \item Conversely, if $g\in\mathcal H$ and $\P_1,\P_2$ are two primitive centered
regular lattice  polytopes such that $g(\P_1)=\P_2$, then $g\in \GL(\Lambda)$ 
 \end{enumerate}
\end{prop}

\begin{proof}
The first point is obvious. 
For the second one, note that if $g\in \mathcal H$ then there exists $r\in \mathbb{Q}$ such that
$r\overrightarrow{g}(\Lambda)=\Lambda$. 
Then let $g\in {\mathcal H},\P_1,\P_2$ be such  $g(\P_1)=\P_2$. 
As $\P_1$ and $\P_2$ are centered,  we have $g(O)=O$ {\it i.e. $g=\overrightarrow{g}$}.
So, we deduce there exist $r\in \mathbb{Q}$ and $h\in\GL(\Lambda)$ such that
$r.h(\P_1)=\P_2$. 
But as $\P_1$ and $\P_2$ are primitive, so $r=1$.
\end{proof}\\

>From now on, we want to classify the primitive centered regular lattice polytopes up to the action 
of $\GL(\L)$.

\section{Root systems}

For root systems, we will use the notation of \cite{Bou}.
Let $\P$ be a \rlp  in $\Lambda_\RR$.
For each edge $a$ of $\P$ with vertices $s_1$ and $s_2$ we consider the subgroup $\RR.\overrightarrow{s_1s_2}\cap\L$
of $\L$ and its two generators $\pm u_a$.
When $a$ runs over all the edges of $\P$, the $\pm u_a$ form a finite subset $\Phi(\P)$ of $\L$.

\begin{prop}
\label{prop:rootsyst}
  The subset $\Phi(\P)$ of $\Lambda_\RR$ is a reduced root system.
\end{prop}

\begin{proof}
  It is clear that $\Phi(\P)$ is finite, does not contain zero, spans $\Lambda_\RR$
and $\ZZ.\alpha\cap\Phi(\P)=\{\pm\alpha\}$ for any $\alpha\in\Phi(\P)$.

Let $\alpha\in\Phi(\P)$ and two vertices $s_1$ and $s_2$ on an edge $a$ parallel to $\alpha$.
Consider a complete flag $\D_1$ of $\P$ starting with $s_1$ and $a$. 
Let $\D_2$ be the complete flag of $\P$ with the same faces from $\D_1$ 
except for the vertex which is $s_2$. 
Let $\sigma\in\Isom(\P)$ such that $\sigma(\D_1)=\D_2$. It is clear that $\overrightarrow{\sigma}$ is a reflection
which maps $\Phi(\P)$ in $\Phi(\P)$ and $\alpha$ on $-\alpha$.

Let $\beta$ be another element of $\Phi(\P)$. The vector $\overrightarrow{\sigma}(\beta)-\beta$ is an element of $\L$
proportional to $\alpha$. Since $\alpha$ has been chosen primitive, $\overrightarrow{\sigma}(\beta)-\beta$ is an entire
multiple of $\alpha$.
\end{proof}\\

The root system $\Phi(\P)$ is said to be {\it associated} to $\P$.
We denote by $\Lambda_P$ and by $\Lambda_R$, respectively the weight and root lattices of 
$\Phi(\P)$. We have:

\begin{prop}
\label{prop:lattice}
  The lattices $\Lambda$, $\Lambda_R$ and $\L_P$  satisfy: $\Lambda_R\subset\Lambda\subset\Lambda_P$.
\end{prop}

\begin{proof}
  The inclusion $\Lambda_R\subset\Lambda$ is obvious.
Let $\alpha\in\Phi(\P)$ and $\lambda\in\Lambda$. We have to prove that $<\lambda,\alpha^\vee>$ 
belongs to $\ZZ$, where $\alpha^\vee$ is the coroot associated to $\alpha$.
But, $\sigma_\alpha(\lambda)=\lambda-<\lambda,\alpha^\vee>\alpha$ belongs to $\Lambda$. 
We can conclude since $\alpha$ is primitive on $\Lambda$.
\end{proof}\\

In the two following propositions, $\P$ is assumed to be centered in $O$.
In this case, $\Isom(\P)$ is a subgroup of $\GL(\Lambda)$.
Let $\Aut(\Phi(\P))=\{g\in\GL(\Lambda_\RR)\ |\ g.\Phi(\P)=\Phi(\P)\}$ denote the automorphism group of $\Phi(\P)$ and $W$ denote the Weyl group of $\Phi(\P)$. Note that $\Aut(\Phi(\P))$ is the semidirect product 
of $W$ and the automorphisms of the Dynkin diagram of $\Phi(\P)$.

\begin{prop}
  \label{prop:isom}
\label{prop:rootsystemsl}
Let $\P$ be a centered regular lattice polytope. We have:
\begin{enumerate}
\item $W\subset\Isom(\P)\subset\Aut(\Phi(\P))$.
\item The lattice $\Lambda$ is stable by $\Isom(\P)$.
\item The root system $\Phi(\P)$ is homogeneous under $\Isom(\P)$.
\end{enumerate}
\end{prop}

\begin{proof}
  The inclusion $W\subset\Isom(\P)$ is a direct consequence of the proof 
of Proposition~\ref{prop:rootsyst}. The rest of the proposition is obvious.
\end{proof}\\

Obviously, $\Isom(\P)$ acts transitively on the set of vertices of $\P$. The following proposition shows a little bit more:

\begin{prop}
\label{prop:Wtrans}
The \rlp $\P$ is assumed to be centered.
The Weyl group $W$ acts transitively on the set of vertices of $\P$.  
\end{prop}

\begin{proof}
Since any edge of $\P$ is parallel to a root of $\Phi(\P)$, any maximal cone of the 
dual fan of $\P$ is an union of Weyl chambers. But, $W$ acts transitively 
on the set of Weyl chambers. The proposition follows.
\end{proof}\\

A face of a regular polytope is a regular polytope. Here, one can say a little bit more:

\begin{prop}
  \label{prop:face}
Let $\F$ be a face of a \rlp $\P$ and $F$ its direction. 
Then, $\F$ is a regular lattice polytope with associated root system $\Phi(\P)\cap F$.
\end{prop}

\begin{proof}
We may assume  that $\P$ is centered.
  It is clear that $\F$ is a lattice regular polytope with root system $\Phi(\F)$ contained in
$\Phi(\P)\cap F$. Let $\alpha\in\Phi(\P)\cap F$: we have to prove that $\alpha$ is parallel to an 
edge of $\F$. 

We claim that the reflection $\sigma_\alpha$ of $W$ associated to $\alpha$ stabilizes $\F$. 
Let $A$ be a point of $\F$. 
The vector $\sigma_\alpha(A)-A$  is collinear to $\alpha$ and so belongs to $F$. 
But, $\F=(A+F)\cap\P$; and so $\sigma_\alpha(A)$ belongs to $\F$.

Consider the kernel $H_\alpha$ of $\sigma_\alpha-\Id$.
Firstly, we assume that $H_\alpha\cap \F$ does not contain any vertex.
Then, there exists an edge $a$ of $\F$ which intersects $H_\alpha$. 
Since $a$ and $\sigma_\alpha(a)$ are edges of $\F$, we have $a=\sigma_\alpha(a)$. 
In particular, $a$ is parallel to $\alpha$.

We now assume that $s$ is a vertex of $\P$ in $H_\alpha\cap \F$. 
Let $b$ be an edge of $\F$ containing $s$. 
Let $\beta$ be a root parallel to $b$. 
This root $\beta$ is neither orthogonal neither collinear to $\alpha$;
so, Proposition~\ref{prop:rootsystemsl} implies that 
$\Phi(\P)\cap{\rm Vect}(\alpha,\beta)$ is a root system of type $A_2$.
Changing $\alpha$ by $-\alpha$ we may assume that $\alpha+\beta$ is a root. 
One easily checks that $\sigma_{\alpha+\beta}(b)$ is parallel to $\alpha$.
\end{proof}\\

\section{Dual Polytope}\label{secdual}

In this section, we define two notions of the dual of a centered regular lattice 
polytope $\P$. Before, we recall the situation in the euclidean case.\\

\subsection{The real case}
 \label{sec:eucliddual}

Let $E$ be a finite dimensional real vector space.
Let $\P$ be a convex polytope in $E$ containing $0$ in its interior.
We denote by $E^*$ the dual of $E$ and set:
$$ 
\P^*=\{\varphi\in E^*{\rm\ s.t.\ }\varphi_{|\P}\geq -1\}.
$$
It is known that $\P^*$ is a convex polytope, called dual of $\P$.
Moreover, $\P^*$ contains $0$ in its interior and the dual $\P^{**}$
of $\P^*$ equals $\P$ modulo the natural identification between $E$ and $E^{**}$.
There is an inclusion-reversing combinatorial correspondence
between the $i$-dimensional faces of $\P$ and the $(n - 1 - i)$-
dimensional faces of $\P^*$.
In particular, if $E$ is euclidean and $\P$ is regular, $\P^*$ is regular too 
with an isomorphic isometry group.\\

Now, we assume that $E$ is euclidean, $\P$ is regular and the barycenter of the vertices 
of $\P$ is $0$.
Consider the convex hull $\P^\vee$ of the barycenters of the facets of $\P$.
With the scalar product, one may identify $E$ and its dual: modulo this identification and 
under our assumptions $\P^*$ are $\P^\vee$ are positively proportional. 
In particular, $\P^\vee$ is regular with the same group as $\P$ and  $\P^{\vee\vee}$ is 
positively proportional to $\P$.\\

The two above constructions of the dual of a regular euclidean polytope can be adapted to
regular lattice polytopes: but the two so obtained notions differ.

\subsection{The lattice  case}

The lattice $\L^*:={\rm Hom}(\L,\ZZ)$ is called the dual of $\L$.
Let  $\P$ be a lattice polytope in $E$ containing $0$ in its interior.
Consider
$$
Q=\{\varphi\in\L^*\otimes\RR{\rm\ s.t.\ }\varphi_{|\P}\geq -1\}.
$$
It is a convex polytope in $\L^*\otimes\RR$ containing $0$ in its interior. 
But, its vertices do not necessarily belong to $\L^*$ but only to $\L^*\otimes\QQ$.
We denote by $\P^*$ the only primitive lattice polytope positively proportional to $Q$.
This lattice polytope $\P^*$ is called the {\it $*$-dual of $\P$}.

Using the properties of $\P^*$ in the real case, one easily check that 
if $\P$ is primitive $\P^{**}=\P$ and that if $\P$ is centered regular so is $\P^*$.\\
 
Now, $\P$ is assumed to be a centered regular lattice polytope. 
There exist a unique positive rational number $k$ such that the barycenters of the vertices 
of the facets of $k.\P$ are primitive vectors in $\L$.
We denote by $\P^\vee$ and call $\vee$-dual of $\P$ the convex hull of its barycenters.

Since $\Isom(\P)$ is finite, there exists a scalar product on $\L_\RR$ such that $\P$ is an 
euclidean regular polytope in $\L_\RR$. 
Then, using the results stated in Section~\ref{sec:eucliddual}, one easily checks that 
if $\P^\vee$ is regular with the same group as $\P$ and that if moreover  
$\P$ is primitive then $\P^{\vee\vee}=\P$.\\

The polytopes $\P^*$ and $\P^\vee$ are not equivalent. 
For example, in dimension two, the two triangles are their
own $\vee$-dual and the $*$-dual one of the other.
In Table~1, we give the $\vee$-dual and $*$-dual of each regular lattice polytope.

\section{Classification}

In this section, we will obtain the classification of the centerd 
regular lattice polytopes.

Let us start by reducing the list of possible root systems.
Let $\P$ be a \cprlp  in $\Lambda_\RR$ of dimension $n$ with associated
root system $\Phi$.
By Proposition~\ref{prop:isom} $\Aut(\Phi)$ acts transitively on  $\Phi$.
Moreover, by Proposition~\ref{prop:face}, there exists a Levi subsystem 
$\Phi'$ of $\Phi$ of rank $n-1$ which is the root system of a regular 
polytope $\Q$ with $\Isom(\Q)$ contained in the stabilizer of $\Phi'$ in 
$\Isom(\P)$.
One easily deduces that the type of $\Phi$ is 
$$A_1^n,\,A_n,\,D_n,\,E_6,\,E_7\  {\rm or\ }E_8.$$

Conversely, let $\Phi$ be a root system in the above list. 
Let us choose a set of simple roots of $\Phi$.
By Proposition~\ref{prop:Wtrans}, the vertices of a centered primitive 
lattice polytope $\P$ with
associated root system $\Phi$ are the orbit by $W$ of a unique dominant 
vertex $s_0$ in $\Lambda_P$. 
Such of polytope $\P$ is also given with a sublattice $\Lambda$ of $\L_P$ 
containing $\L_R$.
Moreover, the polytope $\P$ is completely determined by $\Phi$, $s_0$ 
and $\Lambda$. 
Finally, the polytopes obtained from a pair $(s_0,\,\Lambda)$ and its image
by an automorphism of $\Phi$ are equivalent. 
In Sections~\ref{sec:A1n} to \ref{sec:En}, 
for each possible $\Phi$ we give all the 
possible pairs $(s_0,\Lambda)$ up to the action of $\Aut(\Phi)$.\\

The last step consists to show that each given triple $(\Phi,s_0,\L)$
gives really a regular lattice polytope.
One has to check that the 
stabilizer of $\Lambda$ in $\Aut(\Phi)$ acts transitively on the 
complete flags of the convex hull $\P$ of $W.s_0$ and that $\Phi$ is 
the root system of $\P$. The first verification can be made by checking 
the equality of the cardinality of $\Isom(\P)$ and the set of complete 
flags of $\P$. Using the action of $\Isom(\P)$ the second is equivalent to
check that one root is primitive on one edge of $\P$.
Thereafter, these verifications are left to the reader.

\subsection{Root System $A_1^n$}
\label{sec:A1n}

Here, we assume that the root system $\Phi$ 
associated to the \cprlp $\P$ is of type $A_1^n$.
Let $\om_1,\cdots,\om_n$ be a set of fundamental weights of $\Phi$. 
Then, $\Lambda_P=\ZZ\om_1\oplus\cdots\oplus\ZZ\om_n$ and 
 $\Lambda_R=\ZZ.2\om_1\oplus\cdots\oplus\ZZ.2\om_n$.
Let $k_i\in\ZZ_{>0}$ such that the unique dominant vertex $s_0$ is 
$\sum k_i\om_i$. 
By the action of $W$, the vertices of $\P$ are the $\sum_i \pm k_i\om_i$.
In particular, $\P$ has $n!2^n$ complete flags and $\Isom(\P)=\Aut(\Phi)$.
This implies that  all the $k_i$'s are equal: $s_0=k.\sum_i\om_i$.
Reciprocally, $\Aut(\P)$ acts transitively on the set of flags of 
the convex hull of the $k.\sum_i\pm\om_i$.

Now, we have to determine the possible lattices $\Lambda$.
Necessarily, $\Lambda/\Lambda_R$ is a subgroup of 
$\Lambda_P/\Lambda_R\simeq(\ZZ/2\ZZ)^n$ stable by the action of 
$\Aut(\Phi)$ acting on $(\ZZ/2\ZZ)^n$ by permutations. 
By using for example the canonical bijection between $(\ZZ/2\ZZ)^n$ and the
set of the subsets of $\{1,\cdots,n\}$, one easily checks that the only 
possibilities for $\Lambda$ are:
\begin{enumerate}
\item $\Lambda=\Lambda_R$,
\item $\Lambda=\{\sum_ik_i\om_i\ |\ k_i {\rm \ all\ even\ or\ all\ odd}\}$,
\item $\Lambda=\{\sum_ik_i\om_i\ |\ \sum_ik_i {\rm\ even}\}$,
\item $\Lambda=\Lambda_P$.
\end{enumerate}
But, the edges of the polytopes obtained with $\Lambda=\Lambda_P$ are parallel 
to the $\om_i$'s; so, the root system of $\P$ is $\{\pm\om_i\}$ which is a 
contradiction.
Moreover, for $n=2$, the second and third 
lattices equals. So, we obtain two primitive squares in dimension 2, 
and three primitive cubes for each dimension $n\geq 3$. 

In Table~1, for each choice of $\Lambda$, we give a notation for the class of
the corresponding cube $\C$, the vertex $s_0$, the cardinalities of $\C\cap\L$ 
and of the intersection of $\C$ and an edge of $\L$. We also give the class of
the facets of $\C$ and of its $\vee$ and $*$ duals.  All these results are 
obtained by direct calculations and prove that these cubes are ended 
non equivalent.

\subsection{Root System $D_n$}
\label{sec:Dn}

For convenience, we set $D_1=A_1$, $D_2=A_1\times A_1$ and $D_3=A_3$.

Let $\C^n$ be one of the cubes obtained in the preceding section with $n\geq 2$. 
Its $*$-dual polytope $\CC^n$ is a primitive regular centered lattice polytope 
with $2n$ vertices and $\Isom(\CC^n)$ isomorphic to $\Aut(A_1^n)$.
We easily deduce that the root system of $\CC^n$ is of type $D_n$ and $\Isom(\CC^n)=\Aut(D_n)$, for any
$n\geq 2$.

Since the facets of a cube are  cubes, the stabilizer of $s_0$ in $\Aut(D_n)$ is isomorphic to 
$\Aut(D_{n-1})$; 
then, we may assume that $s_0=k.\om_1$ for a positive integer $k$.

The lattice $\Lambda$ must be stable by the action of $\Aut(D_{n-1})$; 
there are three possibilities (for $n\geq 3$ and with notation of \cite{Bou}):
\begin{enumerate}
\item $\Lambda=\Lambda_R$,
\item $\Lambda=\bigoplus_i\ZZ\varepsilon_i$, 
\item $\Lambda=\Lambda_P$.
\end{enumerate}

So, we obtain three cocubes for each $n\geq 3$ (see Table~1). 

For $n=2$, the two squares are $*$-dual one of the other; in particular,
the cosquares are squares.
\\

Now, let $\P$ be a \cprlp with root system $D_n$ ($n\geq 4$) which is not a 
cocube. Using Proposition~\ref{prop:isom} one easily checks that the root 
system of $\P^\vee$ is $D_n$ too.
The stabilizer of the dominant vertex $s_0$ in $W$ is the stabilizer in $W$ 
of a facet of $\P^\vee$. By Proposition~\ref{prop:face} it is the Weyl group of a 
Levy subsystem of $\Phi$ which is the root system associated to a regular
polytope of dimension $n-1$.
We can deduce that  $n\geq 5$ and $s_0=k.\om_n$ or $s_0=k.\om_{n-1}$, or
$n=4$ and $s_0$ is the multiple of any fundamental weight.

For $n\geq 5$, the two cases  are equivalent using the action of 
$\Aut(\Phi)$. 
We claim that the convex hull $\Q$ of $W.k\om_{n-1}$ is not a regular polytope.
In an adapted base $(e_1,\cdots,e_n)$ of $\Lambda_\RR$, the vertices of $\Q$ are the 
$\sum_i\delta_ie_i$ with $\delta_i=\pm 1$ and 
$\Pi \delta_i=1$.
Let $(x_1,\cdots,x_n)$ denote the dual base of $(e_1,\cdots,e_n)$.
Consider the two linear forms 
$\phi=x_1+\cdots+x_{n-1}-x_{n}$ and $\psi=x_1$. 
The affine hyperplane $\phi=n-2$ define a facet of $\Q$ which is a simplex.
But $\psi=1$ is a facet with $2^{n-2}$ vertices.
So, $\Q$ is not regular since $n\geq 5$.\\

For $n=4$, the three fundamental weights $\om_1$, $\om_3$ and $\om_4$ 
are equivalent modulo $\Aut(D_4)$ and give the cocubes.
Consider the case $s_0=k.\om_2$. 
Since $\om_2$ is the longest root, $k=1$ and the vertices of the convex hull
$\P$ of $W.\om_2$ are the 24 roots of $D_4$ and $\Isom(\P)=\Aut(D_4)$.

Let $(x_1,\cdots,x_4)$ be the dual basis of 
$(\varepsilon_1,\cdots,\varepsilon_4)$ (with notation of \cite{Bou}).
One easily check that the affine hyperplane $\sum x_i=2$ defines 
a facet of $\P$ which is a regular cocube. 
By the action of $\Aut(D_4)$, one obtains the 24 facets:
\begin{itemize}
\item $-2\leq x_1+x_2+x_3+x_4\leq 2$,
\item $-1\leq x_i \leq 1$, for $i=1,\cdots,4$,
\item $-2\leq \sum_{j\neq i} x_j-x_i\leq 2$, for $i=1,\cdots,4$,
\item $x_i+x_j-x_k+x_l\leq 2$, for $\{i,j,k,l\}=\{1,2,3,4\}$ 
and $i<j$ and $k<l$.
\end{itemize}
In particular, $\P$ is regular.

Moreover, $\L_P/\L_R$ is isomorphic to $\ZZ/2\ZZ\times\ZZ/2\ZZ$ and the only lattices 
$\L$ stable by $\Aut(D_4)$ such that $\L_R\subset\L\subset\L_P$ are $\L_P$ and $\L_R$.
So, we obtain two classes of centered primitive regular lattice polytopes called 24-cells 
polytopes.
We denote by $\mathcal{D}_1^4$ those obtained with $\L=\L_P$ and $\mathcal{D}_2^4$ those
obtained with $\L=\L_R$.

The dominant weights in $\P$ are: 
$\om_2$, $\om_1$, $\om_3$, $\om_4$, $\om_1+\om_3$ and $0$. 
By acting $W$ we deduce that $\P\cap\Lambda_P$ contains 
$24+6+6+6+32+1=81$ points and $\P\cap\Lambda_R$ contains $24+1=25$ points.
This gives the cardinality of $\mathcal{D}_i^4\cap\L$, for $i=1$ and $2$.

Since the two 24-cells are the only lattice regular polytopes in dimension four with 
isomorphism group $\Aut(D_4)$ the $\vee$-dual of $\mathcal{D}_1^4$ is either $\mathcal{D}_2^4$
or itself.
But, one easily checks that the barycenter of the facet $\sum_i x_i=2\cap\P$ is 
$\frac12\sum_i\varepsilon_i=\om_4$ and belongs to $\L_P$. 
We deduce that the dual ${\mathcal{D}_1^4}^\vee$ contains strictly less points of $\L$
than $\mathcal{D}_1^4$. Finally, ${\mathcal{D}_1^4}^\vee=\mathcal{D}_2^4$ and 
${\mathcal{D}_2^4}^\vee=\mathcal{D}_1^4$.

\noindent
{\bf Remark}
In this section, we have considered the cocubes for any $n\geq 2$. 
But, the others polytopes with associated root systems of type $D_n$ have only be considered 
when $n\geq 4$. The case $n=2$ has been made in Section~\ref{sec:A1n} and those with $n=3$ will 
be considered in Section~\ref{sec:An}.

\subsection{Root system $A_n$}
\label{sec:simplex}
\label{sec:An}

Consider a \cprlp $\P$ with root system $\Phi$ of type $A_n$ with $n\geq 3$. 
Because of the orders of the Weyl groups, the root system of $\P^\vee$
cannot be of type $E_n$. 
So, we may assume that the root system of $\P^\vee$ is also of type $A_n$; 
if not, we have already meet $\P$.

By Proposition~\ref{prop:face}, the stabilizer of $s_0$ which is the 
stabilizer of a facet of $\P^\vee$ in $\Isom(\P)$ must contains the Weyl 
group of a root system of type $A_{n-1}$ or $A_1\times A_1$ for $n=3$.
This implies that if $n\geq 4$ then $s_0$ equals $k.\om_1$ or $k.\om_n$,
if $n=3$ then $s_0$ is a fundamental weight and implies no restriction 
on $s_0$ if $n=2$.\\

Firstly, we assume that $s_0$ is neither proportional to $\om_1$ or $\om_n$.

Let us fix $n=2$. 
Under our assumption, $\P$ is an hexagon and $\Isom(\P)=\Aut(A_2)$. 
We deduce that $s_0=k.(\om_1+\om_2)$: this gives easily two 
regular hexagons obtained with  $\Lambda$ equal to $\Lambda_R$ and $\Lambda_P$.

If $n=3$, our assumption implies that $s_0=k.\om_2$. 
So, we obtain the three cocubes considered in Section~\ref{sec:Dn}.\\

We now assume that $s_0$ is proportional to $\om_1$. 
The case when $s_0$ is proportional to $\om_n$ is equivalent up to $\Aut(\Phi)$.
The polytope $\P$ is the convex hull of $W.k.\om_1$ that is of 
the $k.\varepsilon_i$'s.
In particular, $\P$ is a simplex and is regular with $\Isom(\P)=W$.

The lattice $\Lambda$ can be any lattice between $\Lambda_R$ and $\Lambda_P$. 
Since $\Lambda_P/\Lambda_R\simeq\ZZ/(n+1)\ZZ$, for each divisor $d$ of $n+1$ we
have exactly one $\Lambda_d$ such that $\Lambda_P/\Lambda_d\simeq\ZZ/d\ZZ$.
For $\Lambda_d$ and $k=d$, $\P$ is a primitive simplex denoted by $\S_d^n$. 
Direct calculation shows that the edges of $\S_d^n$ contain $d+1$ points.
In particular they are pairwise non isomorphic.

The cardinality of $\S_d^n\cap\L$ is a little bit complicated to express. 
For any $\tau\in\ZZ/(n+1)\ZZ$, we denote by $c(\tau)$ the cardinality of the following set
$$
\{
a_1,\cdots,a_{n+1}\in\tau\cap\NN {\rm\ s.t.\ }\sum a_i=d(n+1)
\}.
$$
Then, one easily checks that the cardinality of $\S_d^n\cap\L$ is
$$
\sum_{\tau\in\ZZ/(n+1)\ZZ}c(\tau).
$$
It would be interesting to simplify this formula !

\subsection{Root systems $E_n$}
\label{sec:En}

By absurd, we will prove that there is no \rlp $\P$ with root system
of type $E_6$. We may assume that $\P$ is primitive and centered.
Obviously, the root system of $\P^\vee$ is necessarily $E_6$.
Moreover the root system of a face of $\P$ is either $D_5$ or $A_5$.
The first case is not possible because the regular polytopes with $D_5$ 
as root systems have $\Aut(D_5)$ as isomorphism group. 
But, $\Aut(D_5)$ is not contained in the stabilizer of $D_5$ in $\Aut(E_6)$.
In the second case, one has necessarily $s_0=\om_2$ that is the longest root.
So, the vertices of $\P$ are the roots of $E_6$.
By Proposition~\ref{prop:face}, $\P$ has a facet parallel to the Levy subsystem
of type $A_5$. It follows an easy contradiction.

The same argument shows that there is no regular lattice polytope with 
root system of type $E_7$ and $E_8$.

\section{Description of the regular lattice polytopes}

In the following tabular, for each \cprlp $\P$, we give a notation, its root 
system, the group $\Isom(\P)$ its lattice $\Lambda$, 
its dominant vertex $s_0$, the cardinalities of 
$\P\cap \L$, the number of points in $\L$ on an edge of $\P$, 
\cprlp equivalent to the facets of $\P$ and the duals of $\P^\vee$ and $\P^*$. 
All these elements allow us to distinguish two non isomorphic 
lattice polytopes.

Proofs are given in the preceding section, the others are simple calculation
left to the reader.\\

\begin{remark}
In even dimension more than four, there exist three classes of cocube. 
Two of these three cocubes have the same simplex as facet  and the third another one.
In contradiction in \cite{Kar}, the three cocubes have the same simplex as facet. 
\end{remark}

\bibliographystyle{amsalpha}
\bibliography{biblio}

\newpage
\rotatebox{90}{
\resizebox{\textheight}{0.33\textwidth}{
  \begin{tabular}{c}
\begin{tabular}{|l|c|c|c|c|c|c|c|c|c|c|}
\hline
Type & $\Phi$ & ${\rm Isom}$ & Not.& $\Lambda$ &
$s_0$ & Card. & Edges & Facet & $\P^\vee$ & $\P^*$\\
\hline
\hline
Simplex
&
\begin{tabular}{c}$A_n$\\$n\geq 3$\end{tabular}
&
\begin{tabular}{c}
$W(\Phi)\simeq$\\$\Sym_{n+1}$
\end{tabular}
&
\begin{tabular}{c}
$\mathcal{S}_d^n$\\
for $d | (n+1)$.\\
\end{tabular}
&
\begin{tabular}{c}
$\Lambda_R\subset\Lambda\subset\Lambda_P$\\
with $\#(\Lambda_P/\Lambda)=d$.
\end{tabular}
&
$d\om_1$
&
\begin{tabular}{c}
see\\
Sect \ref{sec:simplex}
\end{tabular}
&
$d+1$&
$\mathcal{S}_n^{n-1}$&
$\mathcal{S}_d^{n}$&
$\mathcal{S}_{\frac{n+1}{d}}^{n}$
\\
\hline
Cubes 
& 
\begin{tabular}{c}
$A_1^n$\\ $n\geq 2$
\end{tabular}  
&
\begin{tabular}{c}
${\rm Aut}(\Phi)$\\$\simeq$\\$(\ZZ/2\ZZ)^n\ltimes\Sym_n$
\end{tabular}
&
\begin{tabular}{c}
$\mathcal{C}_1^n$\\
$\mathcal{C}_2^n$\\
$\left\{
\begin{array}{c}
\mathcal{C}_3^n {\rm\ for\ } n {\rm\ even}\\
~\\
\mathcal{C}_3^n {\rm\ for\ } n{\rm\  odd}
\end{array}
\right.
$
\end{tabular}
&
\begin{tabular}{c}
$\Lambda_R$\\
$k_i\equiv k_j {\rm\ mod\ }2$\\
~\\
$\sum k_i$ even\\
~
\end{tabular}
&
\begin{tabular}{c}
$2\sum\om_i$\\
$\sum\om_i$\\
$\sum\om_i$\\
~\\
$2\sum\om_i$ 
\end{tabular}
&
\begin{tabular}{c}
$3^n$\\
$2^n+1$\\
$\frac{3^n+1}{2}$ \\
~\\
$\frac{5^n-1}{2}$
\end{tabular}
&
\begin{tabular}{c}
$3$\\
$2$\\
$2$ \\
~\\
$3$ 
\end{tabular}
&
\begin{tabular}{c}
$\mathcal{C}_1^{n-1}$\\
$\mathcal{C}_1^{n-1}$\\
~\\
$\mathcal{C}_3^{n-1}$\\
~
\end{tabular}
&
\begin{tabular}{c}
$\mathcal{CC}^n_2$\\
$\mathcal{CC}^n_3$\\
~\\
$\mathcal{CC}^n_1$\\
~
\end{tabular}
&
\begin{tabular}{c}
$\mathcal{CC}^n_2$\\
$\mathcal{CC}^n_3$\\
~\\
$\mathcal{CC}^n_1$\\
~
\end{tabular}
\\
\hline
Cocubes &
\begin{tabular}{c} $D_n$\\ $n\geq 3$\end{tabular} &
\begin{tabular}{c}
${\rm Aut}(\Phi)$\\
$\simeq$\\
$(\ZZ/2\ZZ)^n\ltimes\Sym_n$
\end{tabular}
&
\begin{tabular}{c}
$\mathcal{CC}_1^n$\\
$\mathcal{CC}_2^n$\\
$\left\{
\begin{array}{c}
\mathcal{CC}_3^n {\rm\ for\ } n {\rm\ even}\\
~\\
\mathcal{CC}_3^n {\rm\ for\ } n{\rm\  odd}
\end{array}
\right.$
\end{tabular}
&
\begin{tabular}{c}
$\Lambda_R$\\
$\bigoplus_i\ZZ\varepsilon_i$\\
~\\
$\Lambda_P$\\~
\end{tabular}
&
\begin{tabular}{c}
$2\om_1$\\
$\om_1$\\
~\\$\om_1$\\~
\end{tabular}
&
\begin{tabular}{c}
$4n^2+1$\\
$2n+1$\\
~\\
$2n+1$\\
~
\end{tabular}
&
\begin{tabular}{c}
$3$\\
$2$\\
~\\
$2$\\~
\end{tabular}
&
\begin{tabular}{c}
$\mathcal{S}_n^{n-1}$\\
$\mathcal{S}_n^{n-1}$\\
$\mathcal{S}_{n/2}^{n-1}$\\
~\\
$\mathcal{S}_n^{n-1}$
\end{tabular}&
\begin{tabular}{c}
$\mathcal{C}_3^n$\\
$\mathcal{C}_1^n$\\
~\\
$\mathcal{C}_2^n$\\~
\end{tabular}
&
\begin{tabular}{c}
$\mathcal{C}_3^n$\\
$\mathcal{C}_1^n$\\
~\\
$\mathcal{C}_2^n$\\~
\end{tabular}
\\
\hline
\hline
Hexagon&$A_2$&
\begin{tabular}{c}
${\rm Aut}(\Phi)$\\
$\simeq$\\
$D_6$
\end{tabular}
&
\begin{tabular}{c}
$\mathcal{H}_1^2$\\
$\mathcal{H}_2^2$
\end{tabular}
&
\begin{tabular}{c}
$\Lambda_R$\\
$\Lambda_P$
\end{tabular}
&
$\om_1+\om_2$
&
\begin{tabular}{c}
7\\13
\end{tabular}
&2&&\begin{tabular}{c}
$\mathcal{H}_2^2$\\$\mathcal{H}_1^2$
\end{tabular}&
\begin{tabular}{c}
$\mathcal{H}_1^2$\\$\mathcal{H}_2^2$
\end{tabular}
\\
\hline
24-cell&$D_4$&
\begin{tabular}{c}
${\rm Aut}(\Phi)$\\$\simeq$\\$(\Sym_4\ltimes \ZZ/2\ZZ^3)\ltimes\Sym_3$
\end{tabular}
&
\begin{tabular}{c}
$\mathcal{D}_1^4$\\$\mathcal{D}_2^4$
\end{tabular}
&
\begin{tabular}{c}
$\Lambda_R$\\
$\Lambda_P$
\end{tabular}
&
$\om_2$
&
\begin{tabular}{c}
25\\81
\end{tabular}
&2&
\begin{tabular}{c}
$\CC_1^3$\\$\CC_2^3$
\end{tabular}
&\begin{tabular}{c}
$\mathcal{D}_2^4$\\$\mathcal{D}_1^4$
\end{tabular}&\begin{tabular}{c}
$\mathcal{D}_1^4$\\$\mathcal{D}_2^4$
\end{tabular}\\
\hline
\end{tabular}\\[160pt]    
Remark: We have the following exeptional equalities in dimension two: 
$\mathcal{C}_2^2=\mathcal{C}_3^2$ and 
${\C_1^2}^\vee=C_2^2$.\\[30pt]
 Table 1: List of the centered primitive regular lattice polytopes
\end{tabular}
}
}

\end{document}